\documentclass{amsproc}

\numberwithin{equation}{section}

\usepackage{amssymb}
\usepackage{amsthm}
\usepackage[all]{xy}
\usepackage{graphicx}
\usepackage{textcomp}
\usepackage{color}
\usepackage{dsfont}
\usepackage{slashbox}
\usepackage[T1]{fontenc}

\usepackage{latexsym}
\usepackage{amsmath}
\usepackage{amsfonts}
\usepackage{amssymb}
\usepackage{tikz}
\usepackage[english]{babel}

\newtheorem{conj}{Conjecture}
\newtheorem{lemma}{Lemma}
\newtheorem{prop}{Proposition}
\newtheorem{definition}{Definition}

\newtheorem{remark}{Remark}

\title
[Gog, Magog, and Sch\"utzenberger]
{Gog and  Magog triangles, and the Sch\"utzenberger involution}

\author{Hayat Cheballah}
\address{Laboratoire d'Informatique de Paris Nord, UMR 7030 CNRS, Universit\'e Paris 13,F-93430 Villetaneuse}
\email{Hayat.Cheballah@lipn.univ-paris13.fr}
\author{Philippe Biane}
\address{CNRS, IGM- Universit\'e Paris-Est, 77454 Marne-la-Vall\'e Cedex2, FRANCE}
\email{bian{}e@un{}iv-m{}lv.fr}

\date{\today}

\begin{document}

\begin{abstract}
% \noindent
We describe an approach to finding a bijection between Alternating Sign Matrices and Totally Symmetric Self-Complementary Plane Partitions, which is based on the Sch\"utzenberger involution. In particular we give an explicit
 bijection between Gog and Magog trapezoids with two diagonals.
\end{abstract}

\maketitle

%%%%%%%%%%%%%%%%%%%%%%%%%%%%%%%%%%%%%%%%%%%%%%%%%%%%%%%
\section{Introduction} 

% \noindent
%%%%%%%%%%%%%%%%%%%%%%%%%%%%%%%%%%%%%%%%%%%%%%%%%%%%%%%

\subsection{Alternating Sign Matrices}

An \emph{alternating sign matrix} (ASM) is a square matrix with entries in
$\{-1,0,+1\}$ such that, on each line and on each column, the non zero entries alternate in sign, the sum of each line and
each column being equal to 1. The number of such matrices of size $n$ is
\begin{equation}
\label{eq.AnEnum}
A_n=\prod_{j=0}^{n-1} 
\frac{(3j+1)!}{(n+j)!}
=
1,2,7,42,429,\ldots
\end{equation}
 Zeilberger \cite{zeilberger},  Kuperberg \cite{kuperberg}.
The full story is  in~\cite{bressoud}. 

It has been known for a long time that the numbers $A_n$ also count 
the number of Totally Symmetric Self-Complementary Plane Partitions
(TSSCPP), however no explicit bijection between these classes of objects has been constructed, and finding one is a major open problem in 
combinatorics.

In this paper we propose an approach to this question which is based on the Sch\"utzenberger involution. More precisely, 
we consider Gog and Magog triangles (in the terminology of D. Zeilberger), which are triangular arrays of positive integers, 
satisfying some growth conditions,  in simple bijection  with, respectively, ASMs and TSSCPPs.
The basic idea underlying our approach is that these triangles are examples of Gelfand-Tsetlin patterns to which one can apply 
some known transformations, such as the Sch\"utzenberger involution. In fact we conjecture the existence of a bijection between Gog 
and Magog triangles which can be obtained in two steps, first by making a "modification" of a Gog triangle, based on its inversion pattern, 
then by applying the Sch\"utzenberger involution. This bijection should also preserve trapezoids, which are particular classes of triangles.
 As a first step towards a full bijection we construct here a bijection between $(n,2)$ Gog and Magog trapezoids (the terminology is explained below).

The paper is organized as follows.  

In section 2 we introduce the definitions of Gelfand-Tsetlin triangles, the Gog and Magog triangles and trapezoids, and the Sch\"utzenberger 
involution. Then, in section 3,  we formulate a conjecture on the existence of a bijection between Gog and Magog triangles preserving trapezoids. 
Finally in the last section we give  a bijection between $(n,2)$ Gog and Magog trapezoids.

%%%%%%%%%%%%%%%%%%%%%%%%%%%%%%%%%%%%%%%%%%%%%%%%%%%%%%%
\section{Gog and Magog triangles and trapezoids}
\label{sec.gog_magog}
%%%%%%%%%%%%%%%%%%%%%%%%%%%%%%%%%%%%%%%%%%%%%%%%%%%%%%%
\subsection{Gelfand-Tsetlin}
\begin{definition}
A  Gelfand-Tsetlin  triangle of size $n$ is a triangular array 
$X=(x_{i,j})_{n\geqslant i\geqslant j\geqslant 1}$ of positive integers
$$\begin{matrix}
 x_{n,1}& &x_{n,2}& &\ldots& &x_{n,n-1}& &x_{n,n}\\
  &x_{n-1,1}& &x_{n-1,2}& &\ldots& &x_{n-1,n-1}&  \\
  & &\ldots& &\ldots& &\ldots& & \\
  & & &x_{2,1}& &x_{2,2}& & & \\
  & & & &x_{1,1}& & & & \\
  \end{matrix}
$$
 such that, whenever the  numbers below are defined,
$$x_{i+1,j}\leqslant  x_{i,j}\leqslant x_{i+1,j+1}.$$
\end{definition}
In other words, the triangle is made of $n$  diagonals in the Northwest-Southeast (NW-SE) direction, of lengths 
$n,n-1,\ldots,2,1$ (from  left to right), and it is  weakly increasing in the SE and in the NE directions.

Thus $$\begin{matrix}
 1& &2& &2& &3& &6\\
  &1& &2& &2& &5& \\
  & &2& &2& &4& & \\
  & & &2& &4& & & \\
  & & & &3& & & & \\
 \end{matrix}
$$
is a Gelfand-Tsetlin triangle of size $5$.

Gog and Magog triangles will be obtained from Gelfand-Tsetlin triangles by imposing further conditions on the entries.
\subsection{Gog}

\subsubsection{Triangles} 
\begin{definition} A Gog triangle of size $n$ is a Gelfand-Tsetlin  triangle 
such that

\medskip

$(i)\qquad \qquad \qquad \qquad \qquad  x_{i,j}<x_{i,j+1}, \qquad j<i\leqslant n-1$

\medskip

in other words, such that its rows are strictly increasing, and such that

\medskip

$(ii)\qquad \qquad \qquad \qquad \qquad
x_{n,j}=j,\qquad 1\leqslant j\leqslant n$.

\medskip

\end{definition}

 Here is an example with $n=5$
$$\begin{matrix}
 1& &2& &3& &4& &5\\
  &1& &3& &4& &5& \\
  & &1& &4& &5& & \\
  & & &2& &4& & & \\
  & & & &3& & & & \\
 \end{matrix}
$$

There is  a simple bijection between Gog triangles and Alternating sign matrices (see e.g. \cite{bressoud}). If $(M_{ij})_{1\leqslant i,j\leqslant n}$
is an ASM of size $n$,
then the matrix $\tilde M_{ij}= \sum_{k=i}^nM_{ij}$ has exactly $i-1$ entries 0 and $n-i+1$ entries 1 on row $i$. Let
$(x_{ij})_{j=1,\ldots,i}$ be the columns (in increasing order) with a 1 entry of $\tilde M$ on row $n-i+1$. The triangle 
$X=(x_{ij})_{n\geqslant i\geqslant j\geqslant 1}$ is the Gog triangle corresponding to $M$. 

For example, the above Gog triangle corresponds to the following alternating sign matrix
$$
\begin{pmatrix}
0&1&0&0&0\\
0&0&1&0&0\\
1&-1&0&0&1\\
0&1&-1&1&0\\
0&0&1&0&0
\end{pmatrix}
$$

\subsubsection{Trapezoids}

\begin{definition}
A $(n,k)$ Gog  trapezoid   (for $k\leq n$) is a Gog triangle of size $n$,
$X=(x_{i,j})_{n\geqslant i\geqslant j\geqslant 1}$
such that $x_{i,j}=j$ for $i-j\geq k$.

\end{definition}
Below is a $(5,2)$ Gog trapezoid.
$$\begin{matrix}
1&&2 && 3&& 4&& 5& \\
  &1& &2& &4& &5& & \\
  &&1 & &3& &4& & & \\
  &&&1&&3& & & & \\
&&&&2& & & & 
 \end{matrix}
$$
 See \cite{kratt},  \cite{zeilberger}.
\subsection{Magog}
\subsubsection{Triangles}
\begin{definition}
A Magog triangle of size $n$ is a Gelfand-Tsetlin triangle such that
$$x_{i,i}\leqslant i,\qquad 1\leqslant i\leqslant n.$$

\end{definition}
The set of Magog triangles of size $n$ is in simple bijection with the set of Totally Symmetric Self Complementary Plane Partitions
 (see \cite{bressoud}).

\subsubsection{Trapezoids}
\begin{definition}
A $(n,k)$  Magog trapezoid  (with $k\leqslant n$)  is a Magog triangle
$X=(x_{i,j})_{n\geqslant i\geqslant j\geqslant 1,}$
such that $x_{i,j}=1 $ for $i-j\geq k$.

\end{definition}
 
\subsection{Sch\"utzenberger involution}
\subsubsection{}
Gelfand-Tsetlin triangles label bases of irreducible representations of general linear groups. As such, they are in simple bijection with 
semi-standard Young tableaux. It follows that the Sch\"utzenberger involution, which is defined on SSYTs, can be transported to 
Gelfand-Tsetlin triangles. The following description of this involution has been studied by Berenstein and Kirillov \cite{kirillov}.

First define operators $s_k$, for $k\leqslant n-1$, acting on the set of Gelfand-Tsetlin triangles of size $n$.
If $X=(x_{i,j})_{n\geqslant i\geqslant j\geqslant 1}$ is such a triangle the action of $s_k$ on $X$ is given by 
$s_kX=(\tilde x_{i,j})_{n\geqslant i\geqslant j\geqslant 1}$ with
$$\begin{array}{rcl}
\tilde x_{i,j}&=&x_{i,j},\qquad \text{if\ }i\ne k\\
\tilde x_{k,j}&=&\max(x_{k+1,j},x_{k-1,j-1})+\min(x_{k+1,j+1},x_{k-1,j})-x_{i,j}
\end{array}$$
It is understood that 
$\max(a,b)=\max(b,a)=a$ and $\min(a,b)=\min(b,a)=a$ if the entry $b$ of the triangle is not defined.
The geometric meaning of the transformation of an entry is the following: on row $k$, any entry
$x_{k,j}$ is surrounded by four (or less if it is on the boundary) numbers, increasing from left to right.
$$\begin{matrix}
 x_{k+1,j}& &x_{k+1,j+1}\\
  &x_{k,j}&\\
  x_{k-1,j-1}& &x_{k-1,j} \\
 \end{matrix}
$$
 These four numbers determine a smallest interval containing
$x_{k,j}$, namely $$[\max(x_{k+1,j},x_{k-1,j-1}),\min(x_{k+1,j+1},x_{k-1,j})]$$ and the transformation maps $x_{k,j}$ to its symmetric 
with respect to the center of this interval.

Define $\omega_j=s_js_{j-1}\ldots s_2s_1$.

\begin{definition}

The Sch\"utzenberger involution, acting on Gelfand-Tsetlin triangles of size $n$, is  given by the formula
$$S=\omega_1\omega_2\ldots\omega_{n-1}$$
\end{definition}
It is a non trivial result that $S$ is an involution \cite{kirillov}; beware that the $s_k$ {\sl do not} satisfy the braid relations.

\subsubsection{}
One can  compute the rightmost diagonal of $SX$.
\begin{lemma}\label{SchutzMag}
Let
$X=(X_{i,j})$ be a Gelfand-Tsetlin triangle and $Y=SX$  its image by the Sch\"utzenberger involution, then 

\begin{eqnarray}\label{lonschutz}
 &&Y_{nn}=X_{nn}
\\
&&Y_{kk}=\label{lonschutz2}
\\&&\max_{n=j_0> j_1> j_2\ldots> j_{n-k}\geq 1}\left[
\left(\sum_{i=0}^{n-k-1}X_{j_i+i, j_i}-X_{j_{i+1}+i,j_{i+1}}\right)+X_{j_{n-k}+n-k,j_{n-k}}\right] \nonumber\\
&& \text{for}\ k<n\nonumber
\end{eqnarray}

\end{lemma}
%%%%%%%%%%%%%%%%%%%%%%%%%%%%%%%%%%%%%%%%%%%%%%%%%%%%%%%
{\bf Proof}
We recall the description of the Sch\"utzenberger involution in terms of words and the Robinson-Schensted correspondance.
To the Gelfand-Tsetlin triangle $X$ let us associate the semi-standard Young tableau, with entries in $[1,n]$, such that the shape of the 
tableau formed with letters $u\leq i$ is
the partition $X_{ij},j=1\ldots i$. For example, our Gelfand-Tsetlin triangle
$$\begin{matrix}
 1& &2& &2& &3& &6\\
  &1& &2& &2& &5& \\
  & &2& &2& &4& & \\
  & & &2& &4& & & \\
  & & & &3& & & & \\
 \end{matrix}
$$
corresponds to the tableau (in French notation)

$$\begin{matrix}
 5& && && \\
  4&5& && & \\
  3&3 && && \\
  2& 2& 5&& & \\
  1&1 &1 &2 &4&5 \\
 \end{matrix}
$$
To such a tableau we associate the word $w$ obtained by reading the tableau from top to bottom and from left to right, thus
$$w=5\,4\,5\,3\,3\,2\,2\,5\,1\,1\,1\,2\,4\,5$$ in our example. Then we perform the 
Sch\"utzenberger involution on the word: we read it backwards and replace each letter $i$ by 
$n+1-i$ to give a word $Sw$, in our example
$$Sw=1\,2\,4\,5\,5\,5\,1\,4\,4\,3\,3\,1\,2\,1.$$ 
Observe that this word is a concatenation of nondecreasing words $(Sw)_1,(Sw)_2,\ldots$ corresponding to the successive rows, and 
that these nondecreasing words, viewed as  partitions, are the partitions conjugate to the successive SW-NE diagonals of the original 
Gelfand-Tsetlin triangle.

Applying the Robinson-Schensted algorithm on the word $Sw$ yields an insertion tableau which is  the image of our tableau by the 
Sch\"utzenberger involution.
It is well known that the longest part of the tableau thus obtained is equal to the longest nondecreasing subsequence of the word. 
Thus the largest element of the top row is unchanged.
Moreover, the largest element of the $i^{th}$ row (from bottom) is 
equal to the length of the longest nondecreasing subsequence of the subword of $Sw$ made of numbers $\leq i$.
Now $X_{nn}-X_{j_1j_1}$ is the length of the part of the first nondecreasing subword $(Sw)_1$ of $Sw$ made of letters 
$\leq n-{j_1}$, then $X_{j_1+1,j_1}-X_{j_2+1,j_2}$ is the length of the part of the second nondecreasing subword $(Sw)_2$ of 
$Sw$ made of letters in $[n-j_1,n-j_2]$, and so on.
This yields 
 formula (\ref{lonschutz2}). \qed

\subsubsection{GOGAm triangles}
Since the Sch\"utzenberger involution consists in reading a word backwards and inverting the letters, we introduce the following definition.
\begin{definition}
A GOGAm triangle of size $n$ is a Gelfand-Tsetlin triangle such that its image by the Sch\"utzenberger involution is a Magog triangle of size $n$.
\end{definition}
Thanks to 
Lemma \ref{SchutzMag} we can give a description of GOGAm triangles.
\begin{prop}\label{Schutzgog}
$X=(X_{i,j})$ be a Gelfand-Tsetlin triangle 
then $X$ is a  GOGAm triangle if and only if 
$X_{nn}\leq n$ and, for all $1\leq k\leq n-1$,
and all $n=j_0> j_1> j_2\ldots>j_{n-k}\geq 1$, one has

$$
\left(\sum_{i=0}^{n-k-1}X_{j_i+i, j_i}-X_{j_{i+1}+i,j_{i+1}}\right)+X_{j_{n-k}+n-k,j_{n-k}}\leq k
$$

\end{prop}

{\bf Proof.} Immediate from Lemma \ref{SchutzMag}.
\subsubsection{GOGAm trapezoids}
\begin{definition}
A $(n,k)$ GOGAm trapezoid is a GOGAm  triangle of size $n$  such that
$x_{i,j}=1 $ for $i-j\geq k$. Equivalently, 
 it is the image by the Sch\"utzenberger involution of a $(n,k)$ Magog trapezoid.
\end{definition}
\subsection{A conjecture}

We are now in position to state our conjecture on the Gog-Magog bijection.
\begin{conj}
  There exists a bijection from Gog triangles of size $n$ to Magog triangles of size $n$, which 
  maps $(n,k)$ Gog trapezoids to $(n,k)$ Magog trapezoids for all $k\leq n$. 

\end{conj}

A similar conjecture has been made by Krattenthaler \cite{kratt}.
One can refine the conjecture by considering statistics on Gog and Magog triangles. Such a study will be made in \cite{Cheb}, 
where some further motivation for considering  the Sch\"utzenberger involution will be provided. We consider one such statistic 
in section \ref{stats} below.

In order to construct such a bijection, it is enough to construct a bijection between Gog triangles and GOGAm triangles of 
the same size. In the next section we will give a bijection between $(n,2)$ Gog trapezoids  and   $(n,2)$ GOGAm trapezoids, 
which restricts to  a bijection between $(n,1)$ Gog trapezoids and $(n,1)$ GOGAm trapezoids.

\section{$(n,2)$ Gog and Magog trapezoids}

\subsection{Inversions}  

\begin{definition}

An  inversion in a Gog triangle is a pair $(i,j)$ such that\\ $x_{i,j}=x_{i+1,j}$.

\end{definition}
For example the following Gog triangle contains three inversions, $(2,2)$, $(3,1)$, $(4,1)$, the respective equalities being in red on this picture.

\bigskip

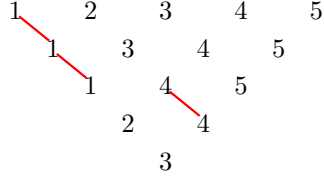
\begin{center}

 \begin{tikzpicture}[scale=0.5]
      \draw  (1,3) node{$1$} ;\draw  (3,3) node{$2$} ;\draw  (5,3) node{$3$} ;\draw  (7,3) node{$4$} ;\draw  (9,3) node{$5$} ;
            \draw  (2,2) node{$1$} ;\draw  (4,2) node{$3$} ;\draw  (6,2) node{$4$} ;\draw  (8,2) node{$5$} ;
                  \draw  (3,1) node{$1$} ;\draw  (5,1) node{$4$} ;\draw  (7,1) node{$5$} ;
                        \draw  (4,0) node{$2$} ;\draw  (6,0) node{$4$} ;
                                 \draw  (5,-1) node{$3$} ;

\draw[line width=0.3mm,color=red] (1.9,2.2)--(1.1,2.85);\draw[line width=0.3mm,color=red] (2.9,1.2)--(2.1,1.85);
\draw[line width=0.3mm,color=red] (5.9,.2)--(5.1,.85);
\end{tikzpicture}

\end{center}
\begin{remark}
The number of inversions of a Gog triangle coincides with the number of inversions of its associated ASM, as defined by Mills, Robbins, Rumsey
\cite{mills}.
\end{remark}

\begin{definition}
Let  $X=(x_{i,j})_{n\geqslant i\geqslant j\geqslant 1}$ be a Gog triangle and let $(i,j)$ be such that
 $1\leqslant i\leqslant j\leqslant n$. 

An inversion $(k,l)$ covers $(i,j)$ if $i=k+p$ and $j=l+p$ for $1\leqslant p\leqslant n-k$. 

\end{definition}
 The entries $(i,j)$ covered by an inversion are depicted with $"+"$ on the following picture.
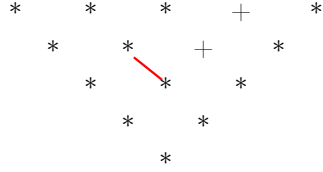
\begin{center}

 \begin{tikzpicture}[scale=0.5]
      \draw  (1,3) node{*} ;\draw  (3,3) node{*} ;\draw  (5,3) node{*} ;\draw  (7,3) node{+} ;\draw  (9,3) node{*} ;
            \draw  (2,2) node{*} ;\draw  (4,2) node{*} ;\draw  (6,2) node{+} ;\draw  (8,2) node{*} ;
                  \draw  (3,1) node{*} ;\draw  (5,1) node{*} ;\draw  (7,1) node{*} ;
                        \draw  (4,0) node{*} ;\draw  (6,0) node{*} ;
                                 \draw  (5,-1) node{*} ;

%\draw[line width=0.3mm,color=red] (2.9,1.2)--(2.1,1.85);
\draw[line width=0.3mm,color=red] (4.15,1.8)--(4.9,1.2);

\end{tikzpicture}

\end{center}

The basic idea for our bijection 
 is that for any inversion in the Gog triangle we should substract 1 to the entries covered by this inversion. This simple minded procedure 
works for $(n,1)$ trapezoids, as we will show as a byproduct of our bijection for   $(n,2)$ trapezoids. It is a good exercise to check this 
directly.   The procedure does not work for $(n,k)$ trapezoids with $k>1$  but, by making some adequate adaptations, we will  obtain  a bijection
 for  trapezoids of size $(n,2)$. 

\subsection{$(n,2)$ trapezoids}
Consider a  $(n,2)$ Gog trapezoid. This  is an array  of the form
$$\begin{array}{ccccccccccccccccc}
 1& &2     &    &  3    &&*   &   &*& &n-2&       &n-1    &      &n\\
  &1&      &2   &      &  *   &    &* &   & n-3 &&b_2    &       &a_1     &  \\
  & &1     &    &    * &       & *  &&  n-4  &&b_3 &       &     a_2  & &  \\
  &       &          &  * & & *&& * &    &*&       & *&&  \\
    &  &    &     &  1    &    &2 &&b_{n-3}&& a_{n-4} &    &&       & & \\
    &    &    &      &      & 1 &  &b_{n-2}&&   a_{n-3}   & &       &&& &  \\
     &   &    &      &       &  &b_{n-1}  &    & a_{n-2}   &&       &    &   & &&  \\
     &  &    &      &        & &   &a_{n-1} &       &&       &      & & && 
\end{array}$$
We shall give an algorithm which builds
a GOGAm triangle from the Gog triangle by successively adding NW-SE diagonals of increasing lengths, and making appropriate changes to the triangle. 
In the end we will obtain a triangle  of the form
$$\begin{array}{ccccccccccccccccc}
 1& &1     &    &  1    &&*   &   &*& &1&       &\beta_1    &      &\alpha_0\\
  &1&      &1   &      &  *   &    &* &   & 1 &&\beta_2    &       &\alpha_1     &  \\
  & &1     &    &    * &       & *  &&  *  &&\beta_3 &       &     \alpha_2  & &  \\
  &       &          &  * & & *&& * &    &*&       & *&&  \\
    &  &    &     &  1    &    &1 &&\beta_{n-3}&& \alpha_{n-4} &    &&       & & \\
    &    &    &      &      & 1 &  &\beta_{n-2}&&   \alpha_{n-3}   & &       &&& &  \\
     &   &    &      &       &  &\beta_{n-1}  &    & \alpha_{n-2}   &&       &    &   & &&  \\
     &  &    &      &        & &   &\alpha_{n-1} &       &&       &      & & && 
\end{array}$$
By Proposition \ref{Schutzgog}, such a triangle is a GOGAm triangle if and only if
\begin{eqnarray*}
\alpha_0&\leq& n\\
 \alpha_0-\alpha_{i}+\beta_{i}&\leq& n-1\quad\text{for}\ 1\leq i\leq n-1\\
\alpha_0-\alpha_i+\beta_i-\beta_j+1&\leq& j-1 \quad\text{for}\ 1\leq i< j\leq n-1
\end{eqnarray*}

\subsection{The algorithm}

First step:
the rightmost NW-SE diagonal consists of one entry $n$ and is not changed, yielding the triangle of size 1 equal to $X^{(1)}=n$.

Second step: The triangle formed by the two first diagonals is
$$\begin{array}{ccc}
 n-1& &n     \\
  &a_1&       \\
\end{array}$$
where $a_1=n\ \text{or} \ n-1$. In the first case, the algorithm yields the triangle 
$$X^{(2)}=\begin{array}{ccc}
 n-1& &n     \\
  &n&       \\
\end{array}$$
in the second case we have an inversion and accordingly substract 1 from the upper right entry, which gives the triangle

$$X^{(2)}=\begin{array}{ccc}
 n-1& &n-1     \\
  &n-1&       \\
\end{array}$$

 Assume now that the first $k$ diagonals have been treated and a triangle $X^{(k)}$ of size $k$, of the form
$$\begin{array}{ccccccccccccc}
 n-k+1& &n-k+1     &            & *&&n-k+1&       &v_1    &      &u_0\\
  &n-k+1&      &*           &     &*&&v_2    &       &u_1     &  \\
   
  &             & *    &  &   * &&   * &          & *&  &\\
  
        &  &        &n-k+1&   &v_{k-2}&   &u_{k-3}    &      &&&  &   \\
        &        &      &   &v_{k-1}  &    & u_{k-2}   &    &  &&  &&  \\
        &          &      &   &   &u_{k-1} &       &       & &    &&  &  
\end{array}$$
has been obtained. 
Furthermore assume that this triangle satisfies the inequalities

\begin{eqnarray}
\label{ineqschutz1}
u_0&\leq& n\\\label{ineqschutz2}
 u_0-u_{i}+v_{i}&\leq& n-1\quad\text{for}\ 1\leq i\leq k-1\\
\label{ineqschutz3} u_0-u_i+v_i-v_j+1&\leq &j-1 \qquad\text{for}\ 1\leq i< j\leq k-1
\end{eqnarray}

and that 
\begin{equation}
u_{k-1}=a_{k-1}.\label{rstat}
\end{equation}

Let us add, on the left of this triangle, the diagonal 
$$\begin{array}{cccccc}
 n-k& &     &    &     &\\
  &n-k&      &   &      &      \\
  & &*     &    &     &       \\
  & &      &n-k   &      &    \\
  & &      &    & v_{k}     &     \\
  & &      &    &      & u_{k}     
\end{array}$$
with $u_k=a_k,\ v_k=b_k$. This yields a triangle $Z^{(k)}$ of size $k+1$ (this triangle will not, in general, be a Gelfand-Tsetlin triangle, 
because the inequality $v_k\leq v_{k-1}$ may be broken).
The algorithm will modify $Z^{(k)}$ to get a triangle $X^{(k+1)}$ of size $k+1$, 
 of the form
$$\begin{array}{ccccccccccccc}
 n-k& &n-k     &       &*&     & *&&n-k&       &v'_1    &      &u'_0\\
  &n-k&      &*           &  &   * &   &*&&v'_2    &       &u'_1     &  \\
  & &*     &    &    * &      & *     &          &  *  &   &u'_2 &&  \\
  & &      &*   &      & *    &  &   * &&   * &          & &  \\
  
  & &      &          &n-k&   &v'_{k-1}&   &u'_{k-2}    &      &&  &   \\
        &    &      &      &   &v'_{k}  &    & u'_{k-1}   &      &&  &&  \\
        &    &      &      &   &   &u'_{k} &       &       &     &&  &  
\end{array}$$

We will  check that the new triangle is a Gelfand-Tsetlin triangle and that 
(\ref{ineqschutz1}), \ldots, (\ref{rstat}) are verified for this new triangle.
The modification will depend  on the inversion pattern in the leftmost diagonal that we have added. In all cases, we will have 
\begin{equation}\label{bottom}
u'_k=u_k
\end{equation}
 the remaining entries being modified as follows, according to the 
 four possibilities for the inversions in the two bottom rows.
\bigskip

$(i)$ The first  case is $v_{k}=n-k,u_k=n-k$, when there are two inversions. Then the modification consists in substracting $1$ 
from each of the entries of the previous triangle, that is we put
$u'_i=u_i-1$, $v'_i=v_i-1$, for $i\leq k-1$, and $v'_{k}=v_k=n-k$.

\bigskip

$(ii)$ The second is the case $v_{k}=n-k<u_k$. Then we  put
$u'_i=u_i$, $v'_i=v_i-1$, for $i\leq k-1$, and $v'_{k}=v_k=n-k$.

\bigskip

$(iii)$
The third case is when $n-k<v_{k}=u_k$. We  put
$u'_i=u_i-1$ for $i\leq k-1$. Observe that $v_{k}=b_{k}=u_k<a_{k-1}=u_{k-1}$, 
therefore $u_i,0\leq i\leq k$ is nonincreasing.  Two cases occur:

\hskip 1cm   $(iiia)$ if the triangle we obtain is a Gelfand-Tsetlin triangle, then  we keep it as the modified triangle, i.e. we put
$v'_i=v_i$ for $i\leq k$. 

\medskip

\hskip 1cm   $(iiib)$ if the triangle is not Gelfand-Tsetlin, then  
  there must exist $j\leq k-1$ with    $v_{j}=u_j$. In this case, we put
$v'_i=v_i-1$, for $i\leq k-1$, and we put
$v'_{k}=n-k$.

\bigskip

$(iv)$ Finally the last case is when $n-k<v_{k}<u_k$. There are two possibilities.

\medskip

\hskip 1cm $(iva)$ if $v_{k}\leq v_{k-1}$, then $Z^{(k)}$ is a Gelfand-Tsetlin triangle, and we do not modifiy it, i.e. we put 
$u'_i=u_i, v'_i=v_i$ for all $i\leq k$, thus
$X^{(k+1)}=Z^{(k)}$.

\medskip

\hskip 1cm $(ivb)$ The last case is  $v_{k}>v_{k-1}$. First we put $u'_i=u_i$ for all $i$.
Let 
\begin{equation}\label{Ldef}
l=\max\{i|v_{k-i}\leq v_k-i\}
\end{equation}
 Since $v_{k-i}$ is nondecreasing and $v_k-i$ is decreasing, one has $l\geq 1$ and $v_{k-i}\leq v_k-i$ for all $i\leq l$. 
We put
$v'_{k}=v'_{k-1}=\ldots=v'_{k-l+1}=n-k$ and $v'_{k-l}=v_{k}-l$, all the other entries being unchanged: $v'_i=v_i$ for $i<k-l$. 

\begin{remark}
Rules $(i),(ii),(iiia),(iva)$ consist just in substracting 1 from entries covered by the inversions in the SE-NW diagonal which has been added. 
The rules $(iiib)$ and $(ivb)$ are more subtle.
\end{remark}
\subsubsection{Proof of the algorithm, first part}
Let us now check that, in each case, we obtain a Gelfand-Tsetlin triangle $X^{(k)}$ satisfying inequalities (\ref{ineqschutz1}), 
(\ref{ineqschutz2}), (\ref{ineqschutz3})
(the identity (\ref{rstat}) is immediate from (\ref{bottom})).

We start with  rules $(i),(ii),(iiia),(iiib),(iva)$.

\bigskip

$(i)$
Since $a_{k-1}=u_{k-1}>v_k$ and $v_{k-1}\geq n-k+1$, $X^{(k+1)}$ is a Gelfand-Tsetlin triangle.
 For $1\leq i< j\leq k-1$ one has $u'_0-u'_i+v'_i-v'_j=u_0-u_i+v_i-v_j$ hence (\ref{ineqschutz3}) is satisfied for these values. Since 
$$u'_0-u'_i+v'_i-v'_k=u_0-u_i+v_i-1-(n-k)\leq n-1-1-(n-k)=k-2$$ we see that  (\ref{ineqschutz3}) is satisfied for all values.
Since $u'_0=u_0-1\leq n-1$ and $u'_i\geq v'_i$ one has (\ref{ineqschutz2})
  and (\ref{ineqschutz1}).

\bigskip

$(ii)$ Again, $X^{(k+1)}$ is clearly a Gelfand-Tsetlin triangle.
For $1\leq i< j\leq k$ we check (\ref{ineqschutz3}) as above, while
(\ref{ineqschutz1}) is clear, 
finally $u'_0-u'_i+v'_i=u_0-u_i+v_i-1\leq n-2$, and
$u'_0-u'_k+v'_k\leq n-1$ since $-u'_k+v'_k\leq -1$, which gives (\ref{ineqschutz2}).

\bigskip

$(iiia)$ Since $u_i>v_i$, one has $u'_i\geq v'_i$ for $i\leq k$, and the triangle $X^{(k+1)}$ is  Gelfand-Tsetlin triangle.

One has,  successively,\begin{eqnarray*}
u'_0&=& u_0-1\\
 u'_0-u'_i+v'_i&=&u_0-u_i+v_i\qquad i<k\\
u'_0-u'_k+v'_k&=&u_0-1\leq n-1\\
u'_0-u'_i+v'_i-v'_j&=&u_0-u_i+v_i-v_j\qquad i<j<k\\
u'_0-u'_i+v'_i-v'_k&<& n-1-(n-k)=k-1 \quad(\text{since }v'_k>n-k)
\end{eqnarray*}
from which inequalities (\ref{ineqschutz1}), (\ref{ineqschutz2}), (\ref{ineqschutz3}) follow.

\medskip

$(iiib)$
The new triangle is clearly Gelfand-Tsetlin.
Furthermore, one has
\begin{eqnarray*}
u'_0&= &u_0-1\\
 u'_0-u'_i+v'_i&=&u_0-u_i+v_i-1\qquad i<k\\
u'_0-u'_k+v'_k&<& u'_0\leq n\\
u'_0-u'_i+v'_i-v'_j&=&u_0-u_i+v_i-v_j\qquad i<j<k\\
u'_0-u'_i+v'_i-v'_k&=& u_0-u_i+v_i-1-(n-k)\leq k-2
\end{eqnarray*}
 which imply inequalities
(\ref{ineqschutz1}), (\ref{ineqschutz2}), (\ref{ineqschutz3}).

\bigskip

$(iva)$ The fact that $X^{(k+1)}$ is Gelfand-Tsetlin is immediate.
 The inequalities are preserved, indeed,  all inequalities involving indices $<k$ are immediate, and one has 
\begin{eqnarray*}
u'_0-u'_k+v'_k&\leq &u'_0-1\leq n-1\quad\text{since }  u'_k>v'_k \\
u'_0-u'_i+v'_i-v'_k&\leq& n-1-(n-k+1)=k-2
\quad \text{since }v'_k>n-k
\end{eqnarray*}

\bigskip
\subsubsection{Proof of the algorithm, second part}
We now consider the last rule,
$(ivb)$.
This is the most delicate part of the proof. We first gather some information on the algorithm which has been constructed up to now.
\begin{lemma}
Just after a step where rule $(i)$ or $(ii)$ is applied, rule $(iiib)$  never  applies.
\end{lemma}
{\it Proof.}
Suppose that rule $(i)$ applies to 
  $Z^{(k)}$, then $n-k=b_k=v_k=a_k=u_k$, and $n-k-1<b_{k+1}=a_{k+1}$ is impossible since this would yield $b_{k+1}\geq a_k$ 
contradicting the Gog strict inequality for the original triangle.
 If rule $(ii)$ applies to $Z^{(k)}$ then 
$v_i<u_i$ in $X^{(k+1)}$, for all $i<k$, therefore rule $(iiib)$ cannot be applied to $Z^{(k+1)}$. 
\qed

\begin{lemma}
If rule $(ivb)$ applies at step $k$, then necessarily at the previous step either rule $(iiib)$ or $(ivb)$ was applied.
\end{lemma} {\it Proof.} If one of the other rules had been applied at the previous step, one would have  $v_{k-1}\geq v_k$.\qed

\begin{lemma}
If rule $(ivb)$ is applied to the triangle $Z^{(k)}$, then
to each of the triangles $Z^{(k-l)},Z^{(k-l+1)}$, $\ldots,Z^{(k-1)}$ either rule   $(iiib)$ or $(ivb)$ was applied.
\end{lemma} {\it Proof.}
Assume that at some step $t<k$ in the algorithm we have  applied rule $(iiia)$ or $(iva)$ to $Z^{(t)}$, then the entry $v_t^{(t+1)}$ 
in the triangle $X^{(t+1)}$ (we emphasize the dependence on the step by adding a superscript) satisfies $b_t=v_t^{(t+1)}$. At   
each next step $s$, we will substract at most   1 from $v_t^{(s)}$, therefore, in the triangle  $Z^{(k)}$, 
 
$$v_t^{(k)}\geq b_t-(k-t-1)\geq b_k-(k-t-1)=v_k^{(k)}-(k-t-1)>v_k^{(k)}+t-k$$ 
It follows that, in $Z^{(k)}$,  one has $l<k-t$ (where $l$ is defined by (\ref{Ldef}). We conclude that, to each of the 
triangles $Z^{(k-l)},Z^{(k-l+1)}$, $\ldots,Z^{(k-1)}$ either rule $(i)$, $(ii)$,  $(iiib)$ or $(ivb)$ was applied. But we have 
seen that rule $(iiib)$ cannot follow  immediately rule $(i)$ or $(ii)$ and that rule $(ivb)$ always follows either rule $(iiib)$ 
or $(ivb)$, so that in fact only rule $(iiib)$ or $(ivb)$ has been applied to each of the triangles $Z^{(k-l)},Z^{(k-l+1)},\ldots,Z^{(k-1)}$.  \qed

\begin{lemma}\label{ivb=}
If rule $(ivb)$ is applied to the triangle $Z^{(k)}$, then
 one has  $$v_{k-1}=\ldots=v_{k-l}=n-k+1.$$ 
\end{lemma} {\it Proof.}
Since $n-k+1\leq v_{k-1}\leq\ldots\leq v_{k-l}$ it is enough to prove that
$v_{k-l}\leq n-k+1$. 
By the preceding Lemma, either rule $(iiib)$ or $(ivb)$ has been applied to the triangles
$Z^{(k-l)},Z^{(k-l+1)}$, $\ldots,Z^{(k-1)}$.
Let us look at the successive values of the entry $v_{k-l}^{(s)}$ in the triangle $X^{(s)}$ (or $Z^{(s)}$).
One has $v_{k-l}^{(k-l+1)}=n-k+l$, since rule $(iiib)$ or $(ivb)$ has been applied to 
$Z^{(k-l)}$.
Each time rule $(iiib)$ is applied   $v_{k-l}^{(s)}$ is decreased by 1. There are two cases

(a) If only rule $(iiib)$ is  applied to $Z^{(k-l)},Z^{(k-l+1)}$, $\ldots,Z^{(k-1)}$ then one has $v_{k-l}^{(k)}=n-k+1$.

(b) If not, let
 $i$ be the least index 
$l\geq i\geq 1$ such that rule $(ivb)$ is applied to $Z^{(k-i)}$, and let 
$l'=\max\{j| v^{(k-i)}_{k-i}-j\geq v^{(k-i)}_{k-i-j}\}$. By rule $(ivb)$ one has 
\begin{equation*}
v^{(k-i+1)}_{k-l'-i}=b_{k-i}-l',\quad v^{(k-i+1)}_{k-i-j}=n-k+i,\ j=0,1,\ldots l'-1
\end{equation*}
Since  rule $(iiib)$ is applied to $Z^{(k-i+1)}$, $\ldots,Z^{(k-1)}$, one has
$
v^{(k)}_{k-l'-i}=b_{k-i}-l'-i+1
$
and
\begin{equation}\label{l'}
 v^{(k)}_{k-p}=n-k+1,\ p=1,2,\ldots l'+i-1
\end{equation}
it follows that 
$$v^{(k)}_{k-l'-i}=b_{k-i}-l'-i+1\geq b_k-l'-i+1=v^{(k)}_{k}-l'-i+1$$ 
hence, by (\ref{Ldef}),
 $$v^{(k)}_{k-l'-i}>v^{(k)}_{k}-l'-i$$ 
 therefore
$l<l'+i$, and $v_{k-l}=n-k+1$ by (\ref{l'}).\qed

  \begin{lemma}
If rule $(iiib)$ or $(ivb)$ is applied to the triangle $Z^{(k)}$, then
 there exists some $i<k-l$ such that $u'_i=v'_i$.
\end{lemma} {\it Proof.}
For rule $(iiib)$ this is easy to see. 

In the case of rule $(ivb)$,
 there exists some step before $k$, when rule $(iiib)$ has been applied and then only rules $(iiib)$ or $(ivb)$ have been applied. 
If rule $(iiib)$ is applied, there must exist an $i$ with $u_i=v_i$, and then applying either rule  $(iiib)$ or $(ivb)$ cannot destroy 
this pair $u_i=v_i$. This implies that there exists some $i$ such that $u'_i=v'_i$. Such a pair cannot exist for $i\geq k-l$ by the preceding 
lemma, therefore $i<k-l$.\qed

\subsubsection{Proof of the algorithm, end}
Assuming that rule $(ivb)$ is applied to the triangle $Z^{(k)}$, 
we can now check that our triangle $X^{(k+1)}$ satisfies all the required properties.
Since $v'_{k-l}=v_k-l$, and $v_{k-l-1}>v_k-l-1$, by the definition of $l$, one has $v'_{k-l-1}\geq v'_{k-l}$. This implies that  $X^{(k+1)}$ 
is  a Gelfand-Tsetlin triangle, as is easily verified.

 Let us check the inequalities
(\ref{ineqschutz1}), (\ref{ineqschutz2}), (\ref{ineqschutz3}).

 First, since $u'_0=u_0$,  (\ref{ineqschutz1}) is clear. Consider $u'_0-u'_i+v'_i$. Since $u'_i=u_i$ is unchanged and $v'_i\leq v_i$   
for all values of $i$, except $v'_{k-l}$, in order to check (\ref{ineqschutz2}) it is enough to consider $u'_0-u'_{k-l}+v'_{k-l}$ and 
$u'_0-u'_k+v'_k$. One has $$u'_{k-l}=u_{k-l}\geq u_k>v_{k}-l=v'_{k-l}$$ therefore
$$u'_0-u'_{k-l}+v'_{k-l}\leq n-1.$$
 Since  $u'_k>v'_k$, one has
$$u'_0-u'_k+v'_k\leq n-1. $$ Consider
$u'_0-u'_i+v'_i-v'_j$, for $i<j\leq k$.

 If $j< k-l$, then  $u'_0-u'_i+v'_i-v'_j=u_0-u_i+v_i-v_j$, so (\ref{ineqschutz3}) is preserved. 

If $j=k-l$, then $u'_i=u_i,\,v'_i=v_i,\,v'_j\geq v_j$ therefore the inequality is again true.

 If $j>k-l>i$, then 
$v'_j=n-k=v_{k-l}-1$ (by Lemma \ref{ivb=}), therefore 
$$u'_0-u'_i+v'_i-v'_j=u_0-u_i+v_i-v_{k-l}+1\leq k-l-1\leq j-2.$$

If $j>k-l=i$ then 
\begin{eqnarray*}
u'_0-u'_i+v'_i-v'_j&=&u_0-u_{k-l}+v_k-l-n+k\\
&=&u_0-n+v_k-u_{k-l}-l+k\leq k-l-1\leq j-2
\end{eqnarray*}
 since $v_k<u_{k-l}$.  

If $k> j>i>k-l$ then
$v'_i-v'_j=v_i-v_j$ and $u'_0-u'_i=u_0-u_i$ therefore
$$u'_0-u'_i+v'_i-v'_j=u_0-u_i+v_i-v_j\leq j-2.$$
Finally
if $k=j>i>k-l$, then
$$u'_0-u'_i+v'_i-v'_k=u_0-u_i+v_i-1-(n-k)\leq n-1-1-(n-k)=k-2.\qed$$ 

\bigskip

Applying the algorithm until we have treated all diagonals, 
we obtain thus a $(n,2)$ GOGAm trapezoid from our $(n,2)$ Gog trapezoid.

\subsubsection{Invertibility}\label{invert}
We can infer from the leftmost SE-NW diagonal of $X^{(k+1)}$ which rule was applied to $Z^{(k)}$. The only ambiguity is whether rule 
$(ii)$, $(iiib)$ or $(ivb)$ has been applied when $n-k=v'_k<u'_k$. Rule $(ii)$ has been applied if and only if one has $u'_i>v'_i$ for all $i<k$. 
 In order to distinguish between rules $(iiib)$ and $(ivb)$  we now state the following lemma.

\begin{lemma}
  Assume $X^{(k+1)}$ is obtained from $Z^{(k)}$ by applying rule $(iiib)$ or $(ivb)$, and let 
$l=1+\max\{i|v'_{k-i}=n-k\}$, then

(a) $v'_{k-l}+l<u'_k$ if rule $(ivb)$ has been applied.

(b) 
$v'_{k-l}+l\geq u'_k$ if rule $(iiib)$ has been applied.
\end{lemma}

Proof. Part (a) is obvious from the statement of rule $(ivb)$, since
$v'_{k-l}+l=v_{k}<u_k=u'_k$.

In order to prove part (b), note that in case $(iiib)$ is applied to $Z^{(k)}$, then by Lemma 2, to all the triangles $Z^{(k-i)}$ for 
$1\leq i\leq l-1$ 
either rule $(iiib)$ or $(ivb)$ has been applied. If only rule $(iiib)$ has been applied to $Z^{(k-l+1)},\ldots,Z^{(k-1)}$, then rule 
 $(iiia)$ or $(iva)$ must have been applied to $Z^{(k-l)}$,
therefore,  
$v'_{k-l}=b_{k-l}-l$ which implies
$v'_{k-l}+l=b_{k-l}\geq b_{k}=a_k=u'_k$.

If rule $(ivb)$ has been applied at some step $t$ with
$k-l+1\leq t\leq k-1$, then let $i$ be the smallest number such that
$(ivb)$ has been applied to $Z^{(k-i)}$. 
By Lemma \ref{ivb=} there exists an $l'\geq 1$ such that
 $$v^{(k-i+1)}_{k-i}=\ldots v^{(k-i+1)}_{k-i-l'+1} =n-k+i-1$$ and 
$$v^{(k-i+1)}_{k-i-l'}=b_{k-i}-l'>n-k+i-1.$$
Since rule $(iiib)$ is applied to $Z^{(k-i+1)},\ldots,Z^{(k-1)}$
it follows that  $$v'_{k}=\ldots = v'_{k-i-l'-1} =n-k$$ and 
$$v'_{k-i-l'}=b_{k-i}-l'-i+1>n-k$$
therefore $l=l'+i$ and   $v'_{k-l}+l\geq u'_k$ since
$v'_{k-l}+l=b_{k-i}\geq b_{k}=a_k=u'_k$.

\qed

\subsection{The inverse map}
\subsubsection{The algorithm}
We now prove that the map defined above has an inverse. 
Let $X$ be a $(n,2)$ GOGAm trapezoid of shape 
$$\begin{array}{ccccccccccccccccc}
 1& &1     &    &  1    &&*   &   &*& &1&       &\beta_1    &      &\alpha_0\\
  &1&      &1   &      &  *   &    &* &   & 1 &&\beta_2    &       &\alpha_1     &  \\
  & &1     &    &    * &       & *  &&  *  &&\beta_3 &       &     \alpha_2  & &  \\
  &       &          &  * & & *&& * &    &*&       & *&&  \\
    &  &    &     &  1    &    &1 &&\beta_{n-3}&& \alpha_{n-4} &    &&       & & \\
    &    &    &      &      & 1 &  &\beta_{n-2}&&   \alpha_{n-3}   & &       &&& &  \\
     &   &    &      &       &  &\beta_{n-1}  &    & \alpha_{n-2}   &&       &    &   & &&  \\
     &  &    &      &        & &   &\alpha_{n-1} &       &&       &      & & && 
\end{array}$$
One has
\begin{eqnarray*}
\alpha_0&\leq& n\\
 \alpha_0-\alpha_{i}+\beta_{i}&\leq& n-1\quad\text{for}\ 1\leq i\leq n-1\\
\alpha_0-\alpha_i+\beta_i-\beta_j+1&\leq& j-1 \quad\text{for}\ 1\leq i< j\leq n-1
\end{eqnarray*}

We shall give an algorithm which is the inverse of the one above.

Let  $k$ be an integer decreasing from $k=n-1$ to $k=0$. Let $Y^{(n)}$ be an empty set, and $X^{(n)}=X$;  at each step we will have a pair 
$(Y^{(k+1)},X^{(k+1)})$ where $Y^{(k+1)}$ 
is an array (non empty only for $k<n-1$)

$$\begin{array}{ccccccccccccccccc}
 1& &2     &    &  *   &&n-k-1  &   && &&       &    &      &\\
  &1&      &2   &      &  *   &    &* &   &  &&   &       &     &  \\
  & &*     &    &    * &       & *  && n-k-1   && &       &      & &  \\
  &       &          &  * & & *&& * &    &b_{k+1}&       & &&  \\
    &  &    &     &  *    &    &2 &&*&& a_{k+1} &    &&       & & \\
    &    &    &      &      & 1 &  &b_{n-2}&&   *   & &       &&& &  \\
     &   &    &      &       &  &b_{n-1}  &    & a_{n-2}   &&       &    &   & &&  \\
     &  &    &      &        & &   &a_{n-1} &       &&       &      & & && 
\end{array}$$
which forms  the leftmost NW-SE diagonals of a Gog triangle, 
and $X^{(k+1)}$ is a Gelfand-Tsetlin triangle: 
$$\begin{array}{ccccccccccccc}
 n-k& &n-k     &            & *&&n-k&       &v'_1    &      &u'_0\\
  &n-k&      &n-k           &     &*&&v'_2    &       &u'_1     &  \\
   
  &             & *    &  &   * &&   * &          & *&  &\\
  
        &  &        &n-k&   &v'_{k-1}&   &u'_{k-2}    &      &&&  &   \\
        &        &      &   &v'_{k}  &    & u'_{k-1}   &    &  &&  &&  \\
        &          &      &   &   &u'_{k} &       &       & &    &&  &  
\end{array}$$
 satisfying the inequalities (\ref{ineqschutz1}), (\ref{ineqschutz2}), (\ref{ineqschutz3}). 
Then we  make a modification of the triangle $X^{(k+1)}$,
according to the rules below, to get a triangle $Z^{(k)}$

$$\begin{array}{ccccccccccccc}
 n-k& &n-k+1     &            & *&&n-k+1&       &v_1    &      &u_0\\
  &n-k&      &n-k+1           &     &*&&v_2    &       &u_1     &  \\
   
  &             & *    &  &   * &&   * &          & *&  &\\
  
        &  &        &n-k&   &v_{k-1}&   &u_{k-2}    &      &&&  &   \\
        &        &      &   &v_{k}  &    & u_{k-1}   &    &  &&  &&  \\
        &          &      &   &   &u_{k} &       &       & &    &&  &  
\end{array}$$
Then we add the leftmost NW-SE  diagonal of this triangle to the right of $Y^{(k+1)}$ to get 
$Y^{(k)}$ (thus $b_{k}=v_{k}$ and $a_{k}=u_{k}$), and take the remaining triangle as $X^{(k)}$.
We will prove that, at each step,   $X^{(k)}$  is a Gelfand-Tsetlin triangle which satisfies the inequalities
(\ref{ineqschutz1}), (\ref{ineqschutz2}), (\ref{ineqschutz3}).   
Also we will prove that, at the next step of the algorithm the entries $a_{k-1},b_{k-1}$ satisfy
\begin{equation}\label{ab's}
n-k+1\leq b_{k-1},\  b_k\leq b_{k-1},\ b_k<a_{k-1},\ b_k\leq a_k\leq a_{k-1}\leq n
\end{equation}
which imply that the triangle $Y^{(0)}$ is a Gog triangle.

We will use the following notation:
if $v'_k=n-k$ and there exists $i<k$ such that $u'_i=v'_i$, then
\begin{equation}\label{el}
l=1+\max\{ j|v'_{k-j}=n-k\}.
\end{equation}
Let us now describe the modification map yilding triangle 
$Z^{(k)}$ from $X^{(k+1)}$ by the inverse algorithm, for which 
we consider several cases, inverse to the cases considered in the direct algorithm.

$(i)$ $n-k=v'_k=u'_k$, then we put
$u_i=u'_i+1,v_i=v'_i+1$ for $i\leq k-1$ and $v_k=v'_k, u_k=u'_k$.

\bigskip

$(ii)$ The second case is $n-k=v'_k<u'_k$, and  $v'_i<u'_i$
for all $i<k$. Then 
we put
$u_i=u'_i, v_i=v'_i+1$, for $i\leq k-1$, and $v_k=v'_k, u_k=u'_k$.
\bigskip

\hskip 1cm $(iiia)$  $n-k<v'_k=u'_k$, then 
we put
$u_i=u'_i+1, v_i=v'_i$ for $i\leq k-1$, and $v_k=v'_k, u_k=u'_k$. 

\medskip

\hskip 1cm $(iiib)$ $n-k=v'_k<u'_k$, there exists $i<k$ such that
$u'_i=v'_i$, and $v'_{k-l}+l\geq u'_k$ (recall (\ref{el})), then we put 
 $u_i=u'_i+1, v_i=v'_i+1$, for $i\leq k-1$, and $v_k=u_k=u'_k$.

\bigskip

\hskip 1cm $(iva)$  $n-k<v'_k<u'_k$, then we put $u_i=u'_i, v_i=v'_i$, $i\leq k$.

\medskip

\hskip 1cm $(ivb)$  $n-k=v'_k<u'_k$, there exists $i<k$ such that
$u'_i=v'_i$, and $v'_{k-l}+l< u'_k$, then we put 
 $u_i=u'_i$,  for $i\leq k$, $v_i=n-k+1$ for $k-l\leq i\leq k-1$,   $v_k=v'_{k-l}+l$, and $v_i=v'_i$ for all other $i$.

\bigskip

Let us now check that this map is well defined. By section
(\ref{invert}), it is an inverse of our modification map.
We consider successively the cases $(i)$,\ldots,$(iv)$  above.
First, by checking all cases successively, one sees that the sequence $a_i$ constructed by the rules above is nonincreasing ($a_i\leq a_{i-1}$), 
and that  $b_i\geq n-i$. The remaining inequalities in (\ref{ab's}) will be checked case by case.
We also have to check that the triangles $X^{(k)}$ are Gelfand-Tsetlin, and that they satisfy (\ref{ineqschutz1}), (\ref{ineqschutz2}), 
(\ref{ineqschutz3}).
The equality (\ref{rstat}) is immediate by inspection.
We start with an observation about rules $(iiib)$ and $(ivb)$.
\begin{lemma}\label{3b4binv}
If rule $(iiib)$ or $(ivb)$ has been applied to the
triangle $X^{(k+1)}$ then in the triangle $X^{(k)}$ there exists a pair
$u_i=v_i$.
\end{lemma}
{\it Proof.} This is immediate for rule $(iiib)$, since adding $1$ to both $u'_i$ and $v'_i$ does not destroy the equality
$u'_i=v'_i$. 

For rule $(ivb)$ we notice that $n-k=v'_k<u'_k$, and
$v'_{k-l}+l<u'_k\leq u'_{k-l}$ imply that $v'_{k-j}<u'_{k-j}$ for $j=1\ldots l$, therefore the inequality $u'_i=v'_i$ must be 
realized for some $i<k-l$, and then $u_i=u'_i=v'_i=v_i$ by rule $(ivb)$. \qed

\bigskip
\subsubsection{Proof of the  algorithm}
We now check successively all rules of the inverse algorithm.

$(i)$ It is clear that the  triangle $X^{(k)}$  is Gelfand-Tsetlin.

We have $u'_0=u'_0-u'_k+v'_k\leq n-1$, 
 this proves (\ref{ineqschutz1}). 

Since 
$u_0-1-u_i+v_i-(n-k)=u'_0-u'_i+v'_i-v'_k\leq k-2$ we have
$u_0-u_i+v_i\leq n-1$. 

All other inequalities
 in (\ref{ineqschutz2}), (\ref{ineqschutz3}) involve differences
like $u_0-u_i$ or $v_i-v_j$ which are not unchanged by the replacement 
$u'\to u, v'\to v$.

Also inequalities (\ref{ab's}) are immediate.

\bigskip

$(ii)$ Since $v'_i<u'_i$ for all $i$ one has $v_i\leq u_i$, hence 
 $X^{(k)}$  is a Gelfand-Tsetlin triangle, and (\ref{ineqschutz1}) is immediate since $u_0=u'_0$. 

Since $u'_0-u'_i+v'_i-v'_k\leq k-2$ one has $u_0-u_i+v_i\leq n-1$, thus 
 (\ref{ineqschutz2}), holds.

Finally  (\ref{ineqschutz3}) comes from $u'_0-u'_i=u_0-u_i$
$v'_i-v'_j=v_i-v_j$.

 The inequalities (\ref{ab's}) at the next step are immediate.

\bigskip

$(iiia)$ Again it is easy to see that $X^{(k)}$ is a Gelfand-Tsetlin triangle.
Since $u'_0-u'_k+v'_k\leq n-1$ and $u'_k=v'_k$ we get $u_0=u'_0+1\leq n$, hence (\ref{ineqschutz1}). 

The other inequalities (\ref{ineqschutz2}), (\ref{ineqschutz3}) are checked similarly.

The inequalities (\ref{ab's})  at the next step are immediate.

\bigskip

$(iiib)$  The fact that $X^{(k)}$ is a Gelfand-Tsetlin triangle 
is immediate.

Since there exists $j$ with $u'_j=v'_j$ one has
$u'_0=u'_0-u'_j+v'_j\leq n-1$ thus $u_0=u'_0+1\leq n$.

Since $u'_0-u'_i+v'_i-v'_k\leq k-2$, it follows that
 $u'_0-u'_i+v'_i\leq n-2$ and 
$u_0-u_i+v_i\leq n-1$. 

The other inequalities are satisfied since
$u'_0-u'_i+v'_i-v'_j=u_0-u_i+v_i-v_j$ for $1\leq i<j\leq k-1$.

We now check the inequalities (\ref{ab's}). 

One has $b_k=a_k<a_{k-1}=u'_{k-1}+1$.

 It remains to see that $b_k\leq b_{k-1}$.

If  $v'_{k-1}>n-k$, then $v_{k-1}=v'_{k-1}+1\geq u'_k=a_k=b_k$ since we are applying rule $(iiib)$ to $X^{(k+1)}$ (in this case, $l=1$). 
At the next step, we will have $b_{k-1}\geq v_{k-1}\geq b_k$.

If  $v'_{k-1}=n-k$, then one has $l>1$, and by Lemma \ref{3b4binv} either rule $(iiib)$ or rule $(ivb)$ applies to $X^{(k)}$.
In either case it is easy to see that $b_k\leq b_{k-1}$.

\bigskip

$(iva)$ In this case, the fact that $X^{(k)}$ is a Gelfand-Tsetlin, as well as the inequalities (\ref{ineqschutz1}),  (\ref{ineqschutz2}), 
(\ref{ineqschutz3}), is immediate. Also the inequalities (\ref{ab's})  are immediate.

\bigskip

$(ivb)$ 
Since $n-k<u'_k\leq u'_{i}$ for $i\leq k-1$ it follows that $u_i\geq n-k+1$ for all $i$.
It is then clear that $X^{(k)}$ is a Gelfand-Tsetlin triangle.

Let us check the  inequalities
(\ref{ineqschutz1}), (\ref{ineqschutz2}), (\ref{ineqschutz3}) for $X^{(k)}$.

Since $u_0=u'_0$  inequality
(\ref{ineqschutz1}) is obvious.

One has  $u_0-u_i+v_i=u'_0-u'_i+v'_i\leq n-1$ for $i< k-l$. For $k>i\geq k-l$ one has $v_i=n-k+1\leq v_{k-l}+l<u_k\leq u_i$ therefore 
$-u_i+v_i\leq -1$ and inequality (\ref{ineqschutz2}) holds.

Inequality $u_0-u_i+v_i-v_j+1=u'_0-u'_i+v'_i-v'_j+1\leq j-1$ holds if $i<j< k-l$. 

If $i<k-l$, one has 
$$
u_0-u_i+v_i-(n-k)+1=u'_0-u'_i+v'_i-v'_{k-l+1}+1\leq k-l$$
hence
$$
 u_0-u_i+v_i-v_{k-l}+1=u_0-u_i+v_i-(n-k+1)+1\leq k-l-1
$$
 which proves (\ref{ineqschutz3}) for $i<j=k-l$.

If $i< k-l< j$ then $v_i=v'_i$ and $v_j\geq v'_j$ therefore the  (\ref{ineqschutz3}) also holds.

 One has
$$u_0-u_{k-l}+v_{k-l}-v_j+1\leq u'_0-u'_{k-l}+v'_{k-l}-v_j+1\leq k-l-1$$ proving  (\ref{ineqschutz3}) for $i=k-l<j$.

If $k-l< i<j$, then $v_i=v_j$, and $v'_i=v'_j$ therefore 
$$u_0-u_i+v_i-v_j+1=u'_0-u'_i+v'_i-v'_j+1\leq j-1.$$

It remain to check inequalities (\ref{ab's}).

After rule $(ivb)$ is applied one has $v_{k-1}=n-k+1$ and, for some $i<k-1$, $u_i=v_i$ therefore rule $(iiib)$ or $(ivb)$ applies 
to the next step. In either case one has $b_{k}<a_{k-1}$.

 Recall that $$b_k=v'_{k-l}+l<u'_k=a_k$$ and 
$$v_{k-1}=\ldots v_{k-l}=n-k+1.$$
 It follows that 
$l'=1+\max\{i|v_{k-1-i}=n-k+1\}\geq l$. 

If $v_{k-1-l'}+l'<u_{k-1}$ then rule $(ivb)$ applies to $X^{(k-1)}$ and $$b_{k-1}=v_{k-1-l'}+l'\geq
v'_{k-l}+l=b_k.$$ 

If $v_{k-1-l'}+l'\geq u_{k-1}$ then $l=l'$,
$u_k=u_{k-1}$ and 
$$b_{k-1}=v_{k-1-l'}+l'=v_{k-1-l}+l=u_{k-1}>b_k.$$
\qed

\subsection{Some properties of the bijection}
\subsubsection{$(n,1)$ trapezoids}
If one starts from a $(n,1)$ trapezoid, then only rules $(i)$ and $(ii)$ apply, and it is easy to see that one gets in the end a 
$(n,1)$ GOGAm trapezoid, and that it is obtained by substracting from any entry of the Gog trapezoid the number of inversions which 
cover it. The same remark applies to the inverse map, so that our bijection restricts to a bijection between $(n,1)$ trapezoids.

\subsubsection{$(n,2,k)$ trapezoids}
A $(n,2,k)$ Gog trapezoid is a $(n,2)$ Gog trapezoid such that
$X_{jj}=n$ for $j\geq k$, and $X_{j,j-1}=\max(X_{k,k-1},j-1)$
for $j\geq k$.

A $(n,2,k)$ Magog trapezoid is a $(n,2)$ Magog trapezoid such that
$X_{ij}=1$ for $j\leq k$.

Our bijection restricts to a bijection between 
$(n,2,k)$ Gog trapezoids and $(n,2,k)$ Magog trapezoids for all $k\leq n$.
\subsubsection{A statistic}\label{stats}
For a Gog triangle $X$ the entry $X_{11}$ gives the position of the 1 in the bottom row of the associated alternating sign matrix.
If $X$ is a $(n,2)$ Gog triangle, it follows from our algorithm  that the $11$ entry of the GOGAm triangle has the same value $X_{11}$. 
From Lemma \ref{SchutzMag}
we conclude that for the $(n,2)$ Magog  triangle $T$, associated
to $X$, one has
$X_{11}=\sum_{i=1}^nT_{i,n}-\sum_{i=1}^{n-1}T_{i,n-1}$. it is known that, more generally, these two statistics on Gog and Magog triangles 
coincide (see e.g. \cite{pzj}).

%%%%%%%%%%%%%%%%%%%%%%%%%%%%%%%%%%%%%%%%%%%%%%%%%%%%%%%

\end{document}